\documentclass{article}
\usepackage{amsfonts}
\usepackage{amsmath}
\usepackage{graphicx}
\usepackage[T1]{fontenc}

\begin{document}

\begin{center}
\bigskip

\vspace*{0.2cm} \thispagestyle{empty}

\bigskip

\bigskip

\bigskip

\bigskip

{\Large B\O URBAKI}

\vspace{0.2cm}

\textit{ou comment Anatole France a r\'{e}invent\'{e} la quantification
existentielle}

\bigskip

\bigskip

\bigskip

\bigskip

{\small par}

\vspace{0.1cm}

F. Zinoun

\vspace{0.3cm}

{\small D\'{e}partement de Math\'{e}matiques}

{\small Facult\'{e} des Sciences de Rabat, Maroc}

{\small f.zinoun@um5r.ac.ma}

{\small Le 12 octobre 2024}

\bigskip

\bigskip

\bigskip

\bigskip

\bigskip

\bigskip

\bigskip 

\bigskip

\bigskip

\includegraphics[scale=0.5]{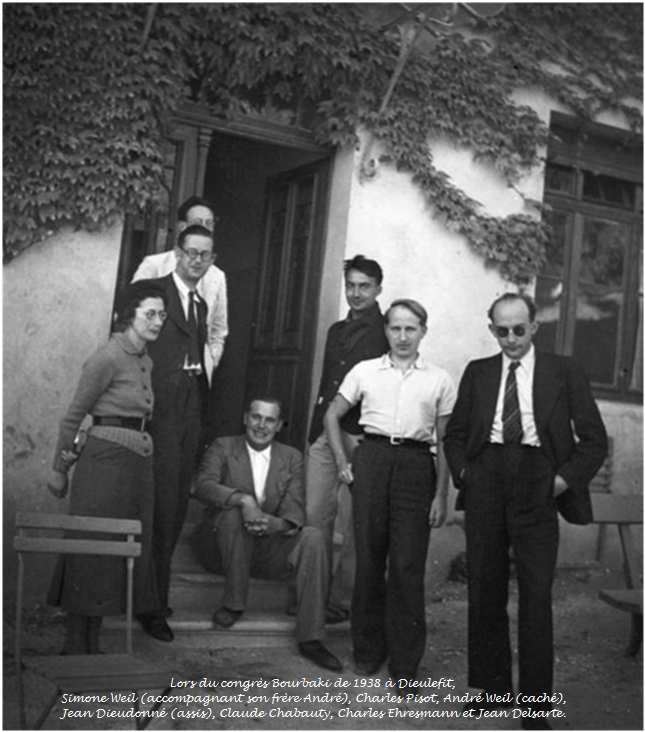}

\newpage

\vspace*{0.2cm}

\bigskip

\bigskip

\bigskip

\bigskip
\end{center}

\textit{Cette histoire est une pure fiction. Toute ressemblance avec des
faits et des personnages existants ou ayant exist\'{e}, ou encore avec une 
\oe uvre litt\'{e}raire connue, ne serait cependant fortuite et ne pourrait 
\^{e}tre le fruit d'une simple co\"{\i}ncidence.}

\begin{center}
\newpage

\bigskip \vspace*{0.2cm}

\bigskip

\bigskip

\bigskip

\bigskip
\end{center}

\textit{\textquotedblleft Bourbaki \'{e}tait. Je puis l'affirmer. Il \'{e}%
tait et il est toujours. Regardez-y et vous vous assurerez qu'\^{e}tre
n'implique nullement la substance et ne signifie que le lien de l'attribut
au sujet (...) Tout le monde croit \`{a} l'existence de Bourbaki. Serais-je
un bon citoyen si je la niais{\tiny \ }? Il faut y regarder \`{a} deux fois avant de
supprimer un article de la foi commune.\textquotedblright }

\newpage

\begin{center}
I

\bigskip

\bigskip
\end{center}

\textit{En visite parisienne, dans le grand appartement des Weil, d'o\`{u}
l'on domine le jardin du Luxembourg, et d'o\`{u} l'on survole la ville
jusqu'au Sacr\'{e}-C\oe ur...}

\bigskip

\guillemotleft {\tiny \ }Ces cahiers de mon enfance, dit M. de Possel, ces petites pages que je garde pr%
\'{e}cieusement depuis Strasbourg, furent pour moi un monde immense, plein
de sourires et de myst\`{e}res.

\vspace{0.1cm}

--- Alain, tu te rappelles Bourbaki{\tiny \ }? demanda Mlle Weil en souriant \`{a} sa
coutume, les l\`{e}vres closes et le nez sur son manuscrit litt\'{e}raire.

\vspace{0.1cm}

--- Bien s\^{u}r Sylvie que je me rappelle Bourbaki{\tiny \ }!... De toutes les
figures qui pass\`{e}rent devant mes yeux quand j'\'{e}tais enfant, celle de
Bourbaki est rest\'{e}e la plus nette dans mon souvenir. Tous les traits de
son visage et de son caract\`{e}re me sont pr\'{e}sents \`{a} la m\'{e}%
moire. Il avait le cr\^{a}ne pointu...

\vspace{0.1cm}

--- Le front bas{\tiny \ }\guillemotright , ajouta Mlle Weil.

\vspace{0.1cm}

Et le fr\`{e}re et la s\oe ur r\'{e}cit\`{e}rent alternativement d'une voix
monotone, avec une gravit\'{e} baroque, les articles d'une sorte de
signalement{\tiny \ }:

\vspace{0.1cm}

\guillemotleft {\tiny \ }Le front bas.

--- Les yeux vairons.

--- Le regard fuyant.

--- Une patte d'oie \`{a} la tempe.

--- Les pommettes aigu\"{e}s, rouges et luisantes.

--- Ses oreilles n'\'{e}taient point ourl\'{e}es.

--- Les traits de son visage \'{e}taient d\'{e}nu\'{e}s de toute expression.

--- Ses mains, toujours en mouvement, trahissaient seules sa pens\'{e}e.

--- Maigre, un peu vo\^{u}t\'{e}, d\'{e}bile en apparence...

--- Il \'{e}tait en r\'{e}alit\'{e} d'une force intellectuelle peu commune.

--- La solution des probl\`{e}mes lui apparaissait en un claquement de
doigts...

--- ... qu'il avait \'{e}normes.

--- Sa voix \'{e}tait tra\^{\i}nante...

--- ... et sa parole mielleuse.{\tiny \ }\guillemotright 

\vspace{0.1cm}

Tout \`{a} coup M. de Possel s'\'{e}cria vivement{\tiny \ }:

\guillemotleft {\tiny \ }Sylvie{\tiny \ }! nous avons oubli\'{e} \textquotedblleft les cheveux jaunes et le
poil rare\textquotedblright .

Recommen\c{c}ons.{\tiny \ }\guillemotright 

\vspace{0.1cm}

Margareth, qui avait entendu avec surprise cette \'{e}trange r\'{e}citation,
demanda \`{a} son p\`{e}re et \`{a} sa tante comment ils avaient pu
apprendre par c\oe ur ce morceau de prose, et pourquoi ils le r\'{e}citaient
comme une litanie.

\vspace{0.1cm}

M. de Possel r\'{e}pondit gravement{\tiny \ }:

\guillemotleft {\tiny \ }Margareth, ce que tu viens d'entendre est un texte consacr\'{e}, je puis
dire liturgique, \`{a} l'usage de la famille de Possel et la famille Weil.
Il convient qu'il soit transmis \`{a} la nouvelle g\'{e}n\'{e}ration, pour
qu'il ne p\'{e}risse pas avec nous. D'ailleurs, j'ai sugg\'{e}r\'{e} \`{a}
ta tante d'inclure ce texte dans son future livre de souvenirs{\tiny \ };
elle juge que, pour entretenir le mythe, il serait plut\^{o}t heureux de le
garder dans les traditions orales de la famille, et j'en conviens. Ta grand-m%
\`{e}re, ma fille, ta grand-m\`{e}re \'{E}velyne, qu'on n'amusait pas avec
des niaiseries, estimait ce morceau, principalement en consid\'{e}ration de
son origine. Elle l'intitula{\tiny \ }: \textit{L'Anatomie de Bourbaki}. Et elle
avait coutume de dire qu'elle pr\'{e}f\'{e}rait, \`{a} certains \'{e}gards,
l'anatomie de Bourbaki \`{a} l'anatomie de Quaresmeprenant.
\textquotedblleft Si la description faite par X\'{e}nomanes, disait-elle,
est plus savante et plus riche en termes rares et pr\'{e}cieux, la
description de Bourbaki l'emporte de beaucoup pour la clart\'{e} des id\'{e}%
es et la limpidit\'{e} du style.\textquotedblright\ Elle en jugeait de la
sorte, sachant qu'elle n'\'{e}tait pas sans conna\^{i}tre l'\OE uvre du docteur Ledouble, o\`{u} il
explique les chapitres trente, trente-un et trente-deux du quart livre de
Rabelais.

\vspace{0.1cm}

--- Je ne comprends pas du tout, dit Margareth.

\vspace{0.1cm}

--- C'est faute de conna\^{\i}tre Bourbaki, ma fille. Il faut que tu saches
que Bourbaki fut la figure la plus famili\`{e}re \`{a} mon enfance et \`{a}
celle de ta tante Sylvie. Chez les Weil, on parlait sans cesse de Bourbaki.
Chacun \`{a} son tour le croyait voir.{\tiny \ }\guillemotright 

\vspace{0.1cm}

Margareth demanda{\tiny \ }:

\guillemotleft {\tiny \ }Qu'est-ce que c'\'{e}tait que Bourbaki{\tiny \ }?{\tiny \ }\guillemotright 

Au lieu de r\'{e}pondre, M. de Possel se mit \`{a} rire, et Mlle Weil aussi
rit, les l\`{e}vres closes.

Margareth portait son regard de l'un \`{a} l'autre. Elle ne trouvait pas 
\'{e}trange que sa tante r\^{\i}t de si bon c\oe ur, mais qu'elle r\^{\i}t
d'accord et en sympathie avec son fr\`{e}re. C'\'{e}tait singulier en effet,
car le demi-fr\`{e}re et la demi-s\oe ur n'avaient pas le m\^{e}me tour
d'esprit.

\vspace{0.1cm}

\guillemotleft {\tiny \ }Papa, dis-moi ce que c'\'{e}tait que Bourbaki. Puisque tu veux que je le
sache, dis-le-moi.

\vspace{0.1cm}

--- Bourbaki, ma fille, \'{e}tait un math\'{e}maticien. Fils d'honorables
scientifiques pold\`{e}ves, il d\'{e}buta charg\'{e} de cours \`{a}
l'Universit\'{e} Royale de Besse-en-Pold\'{e}vie. Mais, fervent d'un
formalisme outrancier, il ne contenta ni ses \'{e}tudiants ni ses pairs.
Ayant quitt\'{e} sa fonction et \'{e}migr\'{e} en Alsace, apr\`{e}s un long p%
\'{e}riple, il y exercera en tant que tuteur priv\'{e} et se consacrera \`{a}
l'\'{e}criture. Ceux qui l'employaient n'eurent cependant pas toujours \`{a}
se louer de lui.{\tiny \ }\guillemotright 

\vspace{0.1cm}

\`{A} ces mots, Mlle Weil, riant encore{\tiny \ }:

\guillemotleft {\tiny \ }Je me rappelle, Alain, que ma m\`{e}re m'avait dit que lorsque mon p\`{e}re ne
trouvait plus sur son bureau une note de l'Acad\'{e}mie, il disait{\tiny \ }:
\textquotedblleft Je soup\c{c}onne Bourbaki d'avoir pass\'{e} par
ici.\textquotedblright\ 

\vspace{0.1cm}

--- Ah{\tiny \ }! dit M. de Possel, Bourbaki n'avait pas une bonne r\'{e}putation.

\vspace{0.1cm}

--- C'est tout{\tiny \ }? demanda Margareth.

\vspace{0.1cm}

--- Non, ma fille, ce n'est pas tout. Bourbaki eut ceci de remarquable,
qu'il nous \'{e}tait connu, familier, et que pourtant...

\vspace{0.1cm}

--- ... il n'existait pas{\tiny \ }\guillemotright, dit Sylvie.

\vspace{0.1cm}

M. de Possel regarda sa s\oe ur d'un air de reproche{\tiny \ }:

\guillemotleft {\tiny \ }Quelle parole, Sylvie{\tiny \ }! et pourquoi rompre ainsi le charme{\tiny \ }? Bourbaki
n'existait pas. L'oses-tu dire, Sylvie{\tiny \ }? Sylvie, le pourrais-tu soutenir{\tiny \ }?
Pour affirmer que Bourbaki n'exista point, que Bourbaki ne fut jamais,
a-t-on assez consid\'{e}r\'{e} les conditions de l'existence et les modes de
l'\^{e}tre{\tiny \ }? Bourbaki existait, ma s\oe ur. Mais il est vrai que c'\'{e}tait
d'une existence particuli\`{e}re.

\vspace{0.1cm}

--- Je comprends de moins en moins, dit Margareth d\'{e}courag\'{e}e.

\vspace{0.1cm}

--- La v\'{e}rit\'{e} t'appara\^{\i}tra clairement tout \`{a} l'heure, ma
fille. Apprends que Bourbaki naquit dans la maturit\'{e} de l'\^{a}ge. J'%
\'{e}tais encore enfant, ta tante devait attendre quelques ann\'{e}es avant
de venir au monde. Mon beau-p\`{e}re Andr\'{e} exer\c{c}ait \`{a} l'Universit%
\'{e} de Strasbourg en tant que professeur de math\'{e}matiques. Ensemble,
avec notre m\`{e}re, ils menaient une vie tranquille et retir\'{e}e dans un
petit appartement sur le bord de l'Ill, jusqu'\`{a} ce qu'ils fussent d\'{e}%
couverts par une vieille dame alsacienne, nomm\'{e}e Mme Shuster, qui vivait
dans son manoir de Brumath, \`{a} quelques lieues de la ville, et qui se
trouva \^{e}tre une grand-tante de notre m\`{e}re. Elle usa d'un droit de
parent\'{e}, mais aussi d'une curieuse passion des math\'{e}matiques, pour
exiger que mon beau-p\`{e}re et ma m\`{e}re vinssent d\^{\i}ner tous les
dimanches chez elle, o\`{u} ils s'ennuyaient excessivement. Elle disait
qu'il \'{e}tait honn\^{e}te de d\^{\i}ner en famille le dimanche et que
seuls les gens mal n\'{e}s n'observaient pas cet ancien usage. Mon beau-p%
\`{e}re pleurait d'ennui et se lassait de ces interminables discussions sur
l'avenir des math\'{e}matiques et de leur enseignement. D\'{e}j\`{a}, il
devait de temps \`{a} autre rentrer \`{a} Paris pour son s\'{e}minaire,
entre autres r\'{e}unions, ce qui ne lui laissait point de force pour
s'infliger d'autres d\'{e}placements. Son d\'{e}sespoir faisait peine \`{a}
voir. Mais Mme Shuster ne le voyait pas. Elle ne voyait rien. Ma m\`{e}re
avait plus de courage. Elle souffrait autant que mon beau-p\`{e}re, et peut-%
\^{e}tre davantage, et elle souriait.

\vspace{0.1cm}

--- Les femmes sont faites pour souffrir, dit Sylvie avec \'{e}motion, en
tournant tendrement les yeux vers le portrait de sa tante, Simone Weil,
qu'on croirait sa s\oe ur jumelle, sous son b\'{e}ret basque.

\vspace{0.1cm}

--- Sylvie, tout ce qui vit au monde est destin\'{e} \`{a} la souffrance. En
vain les Weil refusaient ces funestes invitations. La voiture de Mme Shuster
venait les prendre chaque dimanche, apr\`{e}s midi. Il fallait aller \`{a}
Brumath{\tiny \ }; c'\'{e}tait une obligation \`{a} laquelle il \'{e}tait
absolument interdit de se soustraire. C'\'{e}tait un ordre \'{e}tabli, que
la r\'{e}volte pouvait seule rompre. Mon beau-p\`{e}re enfin se r\'{e}volta,
et jura de ne plus accepter une seule invitation de Mme Shuster, laissant 
\`{a} ma m\`{e}re le soin de trouver \`{a} ces refus des pr\'{e}textes d\'{e}%
cents et des raisons vari\'{e}es, c'est ce dont elle \'{e}tait le moins
capable. Notre m\`{e}re ne savait pas feindre.

\vspace{0.1cm}

--- Dis, Alain, qu'elle ne voulait pas. Elle aurait pu mentir comme les
autres.

\vspace{0.1cm}

--- Il est vrai de dire que lorsqu'elle avait de bonnes raisons, elle les
donnait plut\^{o}t que d'en inventer de mauvaises. Tu sais, ma s\oe ur,
qu'il lui arriva un jour de dire, \`{a} table{\tiny \ }: \textquotedblleft
Heureusement qu'Alain a la coqueluche{\tiny \ }: nous n'irons pas de longtemps \`{a}
Brumath.\textquotedblright\ 

\vspace{0.1cm}

--- C'est pourtant vrai, je crois{\tiny \ }! dit Sylvie.

\vspace{0.1cm}

--- Je gu\'{e}ris, Sylvie. Et Mme Shuster vint dire un jour \`{a} notre m%
\`{e}re{\tiny \ }: \textquotedblleft Ma mignonne, je compte bien que vous viendrez
avec votre mari d\^{\i}ner dimanche \`{a} Brumath.\textquotedblright\ Notre m%
\`{e}re, charg\'{e}e express\'{e}ment par son mari de pr\'{e}senter \`{a}
Mme Shuster un valable motif de refus, imagina, en cette extr\'{e}mit\'{e},
une raison qui n'\'{e}tait pas v\'{e}ritable. \textquotedblleft Je regrette
vivement, ch\`{e}re madame. Mais cela nous sera impossible. Dimanche, nous
attendons un tuteur.\textquotedblright

\vspace{0.1cm}

\textquotedblleft Vous attendez un tuteur{\tiny \ }! Pourquoi{\tiny \ }? --- Pour initier
Alain aux rudiments de l'arithm\'{e}tique.\textquotedblright

\vspace{0.1cm}

\guillemotleft {\tiny \ }Et ma m\`{e}re, ayant tourn\'{e} involontairement les yeux vers son mari,
qu'elle venait de lui substituer un enseignant de math\'{e}matiques,
reconnut avec effroi l'invraisemblance de son invention. \textquotedblleft
M. Weil, dit Mme Shuster, pourrait bien r\'{e}viser avec son beau-fils le
long de la semaine. D'ailleurs, cela vaudrait mieux... --- Alain pr\'{e}f%
\`{e}re que ce soit quelqu'un d'autre, c'est plus d\'{e}contractant pour
lui.\textquotedblright 

\vspace{0.1cm}

\guillemotleft {\tiny \ }J'ai remarqu\'{e} souvent que les raisons les plus absurdes et les plus
saugrenues sont les moins combattues{\tiny \ }: elles d\'{e}concertent l'adversaire.
Mme Shuster insista, moins qu'on ne pouvait l'attendre d'une personne aussi
peu dispos\'{e}e qu'elle \`{a} d\'{e}mordre. En se levant de dessus son
fauteuil, elle demanda{\tiny \ }: \textquotedblleft Comment l'appelez-vous, ma
mignonne, votre prof{\tiny \ }? --- Bourbaki\textquotedblright , r\'{e}pondit ma m%
\`{e}re sans h\'{e}sitation.

\vspace{0.1cm}

\guillemotleft {\tiny \ }Le nom lui vint naturellement \`{a} l'esprit en se rappelant des canulars
qui lui ont \'{e}t\'{e} cont\'{e}s par son mari, le plus spectaculaire \'{e}%
tant celui jou\'{e}, quelques lustres auparavant, par un certain M. Husson,
alors \'{e}tudiant \`{a} l'\'{E}cole Normale Sup\'{e}rieure.
Vraisemblablement emprunt\'{e} \`{a} un g\'{e}n\'{e}ral napol\'{e}onien,
c'est aussi le nom sous lequel son mari, avec une bande de jeunes math\'{e}%
maticiens, entendaient r\'{e}ellement \'{e}diter des textes de math\'{e}%
matique(s) dans les ann\'{e}es \`{a} venir. Mais Bourbaki, le professeur
strasbourgeois, \'{e}tait ce jour-l\`{a} nomm\'{e}. D\`{e}s lors il exista.
Lui sonnant comme un nom slave, ou plut\^{o}t grec, Mme Shuster s'en alla en
ronchonnant : \textquotedblleft Bourbaki{\tiny \ }! Il me semble bien que je connais 
\c{c}a. Bourbaki{\tiny \ }? Bourbaki{\tiny \ }! Je ne connais que lui. Mais je ne me rappelle
pas... O\`{u} exerce-t-il{\tiny \ }? --- Nulle part. Quand on a besoin de lui, on le
lui fait dire chez l'un ou chez l'autre. --- Ah{\tiny \ }! je le pensais bien{\tiny \ }: un
fain\'{e}ant et un vagabond... un rien du tout. M\'{e}fiez-vous de lui et de
ses le\c{c}ons, ma mignonne.\textquotedblright 

\vspace{0.1cm}

\guillemotleft {\tiny \ }Bourbaki avait d\'{e}sormais un caract\`{e}re.{\tiny \ }\guillemotright 

\bigskip

\bigskip

\begin{center}
II

\bigskip

\bigskip
\end{center}

\noindent Nicolette Weil \'{e}tant survenue, M. de Possel la mit au point de
la conversation{\tiny \ }:

\vspace{0.1cm}

\guillemotleft {\tiny \ }Nous parlions de celui qu'un jour notre m\`{e}re fit na\^{\i}tre professeur 
\`{a} Strasbourg et nomma par son nom. D\`{e}s lors il agit.

\vspace{0.1cm}

--- Cher fr\`{e}re, voudrais-tu r\'{e}p\'{e}ter{\tiny \ }? dit Mlle Weil en essuyant
le verre de ses lunettes.

\vspace{0.1cm}

--- Volontiers, r\'{e}pondit M. de Possel. Il n'y avait pas de professeur.
Le professeur n'existait pas. Ma m\`{e}re dit{\tiny \ }: \textquotedblleft Nous
attendons un tuteur\textquotedblright . Aussit\^{o}t le tuteur fut. Et il
agit. M. Husson dit{\tiny \ }: \textquotedblleft J'\'{e}nonce un th\'{e}or\`{e}me de
Bourbaki\textquotedblright . Aussit\^{o}t l'auteur fut. Et il agit.

\vspace{0.1cm}

--- Cher fr\`{e}re, demanda Mlle Weil, comment agit-il, puisqu'il n'existait
pas?

\vspace{0.1cm}

--- Il avait une sorte d'existence, r\'{e}pondit M. de Possel.

\vspace{0.1cm}

--- Tu veux dire une existence imaginaire, Alain, r\'{e}pliqua Mlle Weil.

\vspace{0.1cm}

--- N'est-ce donc rien qu'une existence imaginaire{\tiny \ }? Et les personnages
mythiques ne sont-ils donc pas capables d'agir sur les hommes{\tiny \ }? R\'{e}fl\'{e}%
chissez sur la mythologie mesdemoiselles, et vous vous apercevrez que ce
sont, non point des \^{e}tres r\'{e}els, mais des \^{e}tres imaginaires qui
exercent sur les \^{a}mes l'action la plus profonde et la plus durable.
Partout et toujours des \^{e}tres, qui n'ont pas plus de r\'{e}alit\'{e} que
Bourbaki, ont inspir\'{e} aux peuples la haine et l'amour, la terreur et
l'esp\'{e}rance, conseill\'{e} des crimes, re\c{c}u des offrandes, fait les m%
\oe urs et les lois. R\'{e}fl\'{e}chissez sur l'\'{e}ternelle mythologie.
Bourbaki est un personnage fictif, des plus obscurs, j'en conviens. Le
grossier satyre, assis jadis \`{a} la table de nos paysans du Nord, fut jug%
\'{e} digne de para\^{\i}tre dans un tableau de Jorda\"{e}ns et dans une
fable de La Fontaine. Le fils velu de Sycorax entra dans le monde sublime de
Shakespeare. Bourbaki, moins heureux, sera toujours m\'{e}pris\'{e} des
artistes et des po\`{e}tes. Il lui manque la grandeur et l'\'{e}tranget\'{e}%
, le style et le caract\`{e}re. Il naquit dans des esprits trop
raisonnables, parmi des gens qui, en qu\^{e}te de V\'{e}rit\'{e}, ne
faisaient qu'argumenter et n'avaient point cette imagination charmante qui s%
\`{e}me les fables. Je pense, mesdemoiselles, que j'en ai dit assez pour
vous faire conna\^{\i}tre la v\'{e}ritable nature de Bourbaki.

\vspace{0.1cm}

--- Je la con\c{c}ois{\tiny \ }\guillemotright, dit Nicolette.

\vspace{0.1cm}

Et M. de Possel poursuivit son discours{\tiny \ }:

\guillemotleft {\tiny \ }Bourbaki \'{e}tait. Je puis l'affirmer. Il \'{e}tait et il est toujours.
Regardez-y, mesdemoiselles, et vous vous assurerez qu'\^{e}tre n'implique
nullement la substance et ne signifie que le lien de l'attribut au sujet,
n'exprime qu'une relation.

\vspace{0.1cm}

--- Sans doute, dit Nicolette, mais \^{e}tre sans attributs c'est \^{e}tre
aussi peu que rien. Je ne sais plus qui a dit autrefois{\tiny \ }: \textquotedblleft
Je suis celui qui est.\textquotedblright\ Excusez-moi, on ne peut tout se
rappeler. Mais l'inconnu qui parla de la sorte commit une rare imprudence.
En donnant \`{a} entendre par ce propos inconsid\'{e}r\'{e} qu'il \'{e}tait d%
\'{e}pourvu d'attributs et priv\'{e} de toutes relations, il proclama qu'il
n'existait pas et se supprima lui-m\^{e}me \'{e}tourdiment. Je parie qu'on
n'a plus entendu parler de lui.

\vspace{0.1cm}

--- Tu as perdu, r\'{e}pliqua Alain. Il a corrig\'{e} le mauvais effet de
cette parole \'{e}go\"{\i}ste en s'appliquant des pot\'{e}es d'adjectifs, et
l'on a beaucoup parl\'{e} de lui, le plus souvent sans aucun bon sens.

\vspace{0.1cm}

--- Je ne comprends pas, dit Nicolette.

\vspace{0.1cm}

--- Il n'est pas n\'{e}cessaire de comprendre{\tiny \ }\guillemotright, intervint Margareth.

\vspace{0.1cm}

Et elle pria son p\`{e}re de continuer de lui parler de Bourbaki.

\vspace{0.1cm}

\guillemotleft {\tiny \ }Tu es bien aimable de me le demander, ma ch\'{e}rie, fit le père.

\guillemotleft {\tiny \ }Bourbaki naquit dans la premi\`{e}re moiti\'{e} de notre si\`{e}cle, \`{a}
Strasbourg, ou peut-\^{e}tre \`{a} Paris, rue d'Ulm, ou encore \`{a}
Aligarh, en Inde, on ne sait pas vraiment. Il lui aurait peut-\^{e}tre mieux
valu na\^{\i}tre des si\`{e}cles auparavant dans un lieu des Mille et Une
Nuits. \c{C}'aurait \'{e}t\'{e} alors un malin esprit d'une merveilleuse
adresse.

\vspace{0.1cm}

--- Une tasse de th\'{e}, Nicolette{\tiny \ }? dit Sylvie.

\vspace{0.1cm}

--- Bourbaki \'{e}tait-il donc un malin esprit{\tiny \ }? demanda Nicolette, avec son
accent am\'{e}ricain.

\vspace{0.1cm}

--- Il \'{e}tait malin, r\'{e}pondit M. de Possel, il l'\'{e}tait en quelque
mani\`{e}re, mais il ne l'\'{e}tait pas absolument. Il en est de lui comme
des diables qu'on dit tr\`{e}s m\'{e}chants, mais en qui l'on d\'{e}couvre
de bonnes qualit\'{e}s quand on les fr\'{e}quente. Et je serais dispos\'{e} 
\`{a} croire qu'on a fait tort \`{a} Bourbaki. Mme Shuster qui, pr\'{e}venue
contre lui, l'avait tout de suite soup\c{c}onn\'{e} d'\^{e}tre un fain\'{e}%
ant, un vagabond et un moins que rien, commen\c{c}ait \`{a} montrer une
certaine estime \`{a} son \'{e}gard lorsque ma m\`{e}re avait commis
l'indiscr\'{e}tion de lui montrer un fascicule qui allait devenir le premier
volume des \'{E}l\'{e}ments de math\'{e}matique de Bourbaki. Elle r\'{e}fl%
\'{e}chit aussi que puisque ma m\`{e}re l'employait, elle qui n'\'{e}tait
pas riche, et puisque le bonhomme se consacre \`{a} la r\'{e}daction de
notes et d'ouvrages ne serait-ce que par amour de la discipline, c'\'{e}tait
qu'il se contentait de peu, et elle se demanda si elle n'aurait pas avantage 
\`{a} le faire travailler dans une \'{e}cole \`{a} Brumath dont elle d\'{e}tenait de bons parts, pr\'{e}f\'{e}rablement \`{a} certains enseignants qui
avaient meilleurs renoms, mais aussi plus d'exigences. C'\'{e}tait la rentr%
\'{e}e. Elle pensa que si Mme Weil, qui \'{e}tait pauvre, ne donnait pas
grand-chose \`{a} Bourbaki, elle-m\^{e}me, qui \'{e}tait riche, lui
donnerait moins encore, puisque c'est l'usage que les riches paient moins
cher que les pauvres. Et elle voyait d\'{e}j\`{a} l'\'{e}cole gagner un
auxiliaire sans qu'on y f\^{\i}t grande d\'{e}pense. \textquotedblleft
J'aurai l'\oe il, se dit-elle, \`{a} ce que Bourbaki ne fl\^{a}ne point. Je
ne risque rien et ce sera tout profit. Ces vagabonds se montrent quelquefois
plus d\'{e}vou\'{e}s dans leur t\^{a}che que les
permanents.\textquotedblright\ Elle r\'{e}solut d'en faire l'essai et dit 
\`{a} ma m\`{e}re : \textquotedblleft Mignonne, envoyez-moi Bourbaki. Je le
ferai travailler \`{a} Brumath.\textquotedblright\ Ma m\`{e}re le lui
promit. Elle l'e\^{u}t fait volontiers. Mais vraiment ce n'\'{e}tait pas
possible. Mme Shuster attendit Bourbaki \`{a} Brumath, et l'attendit en
vain. Elle avait de la suite dans les id\'{e}es et de la constance dans ses
projets. Quand elle revit ma m\`{e}re, elle se plaignit \`{a} elle de
n'avoir pas de nouvelles de Bourbaki. \textquotedblleft Mignonne, vous ne
lui avez donc pas dit que je l'attendais{\tiny \ }? --- Si{\tiny \ }! mais il est \'{e}trange,
bizarre... --- {\tiny \ }Oh{\tiny \ }! je connais ce genre-l\`{a}. Je le sais par c\oe ur votre
Bourbaki. Mais il n'y a pas d'enseignant assez lunatique pour refuser de
venir exercer dans l'une des plus belles r\'{e}gions d'Alsace. Notre \'{e}%
cole est connue, je pense. Dites-moi seulement o\`{u} il loge{\tiny \ };
j'irai moi-m\^{e}me le trouver.\textquotedblright\ Ma m\`{e}re r\'{e}pondit
qu'elle ne savait pas o\`{u} logeait Bourbaki, qu'on ne lui connaissait pas
de domicile, qu'il \'{e}tait sans feu ni lieu. \textquotedblleft Je ne l'ai
pas revu, madame. Je crois qu'il se cache.\textquotedblright\ Pouvait-elle
mieux dire{\tiny \ }?

\vspace{0.1cm}

\guillemotleft {\tiny \ }Mme Shuster pourtant ne l'\'{e}couta pas sans d\'{e}fiance{\tiny \ }; elle
la soup\c{c}onna de circonvenir Bourbaki, de le soustraire aux recherches,
dans la crainte de le perdre ou de le rendre plus exigeant. Et elle la jugea
vraiment trop \'{e}go\"{\i}ste. Beaucoup de jugements accept\'{e}s par tout
le monde, et que l'histoire a consacr\'{e}s, sont aussi bien fond\'{e}s que
celui-l\`{a}.

\vspace{0.1cm}

--- C'est pourtant vrai, dit Margareth.

\vspace{0.1cm}

--- Qu'est-ce qui est vrai{\tiny \ }? demanda Sylvie \`{a} demi sommeillant, son manuscrit toujours entre les mains.

\vspace{0.1cm}

--- Que les jugements de l'histoire sont souvent faux. Je me souviens, papa,
que tu as dit un jour{\tiny \ }: \textquotedblleft Mme Roland \'{e}tait bien na\"{\i}%
ve d'en appeler \`{a} l'impartiale post\'{e}rit\'{e} et de ne pas
s'apercevoir que, si ses contemporains \'{e}taient de mauvais singes, leur
post\'{e}rit\'{e} serait aussi compos\'{e}e de mauvais
singes.\textquotedblright

\vspace{0.1cm}

--- Margareth, demanda Sylvie, quel rapport y a-t-il entre l'histoire de
Bourbaki et ce que tu nous contes l\`{a}{\tiny \ }?

\vspace{0.1cm}

--- Un tr\`{e}s grand, ma tante.

\vspace{0.1cm}

--- Je ne le saisis pas.{\tiny \ }\guillemotright 

\vspace{0.1cm}

M. de Possel, qui n'\'{e}tait pas ennemi des digressions, r\'{e}pondit \`{a}
sa fille{\tiny \ }:

\guillemotleft {\tiny \ }Si toutes les injustices \'{e}taient finalement r\'{e}par\'{e}es en ce
monde, on n'en aurait jamais imagin\'{e} un autre pour ces r\'{e}parations.
Comment voulez-vous que la post\'{e}rit\'{e} juge \'{e}quitablement tous les
morts{\tiny \ }? Comment les interroger dans l'ombre o\`{u} ils fuient{\tiny \ }? D\`{e}s
qu'on pourrait \^{e}tre juste envers eux, on les oublie. Mais peut-on jamais 
\^{e}tre juste{\tiny \ }? Et qu'est-ce que la justice{\tiny \ }? Mme Shuster, du moins, fut
bien oblig\'{e}e de reconna\^{\i}tre \`{a} la longue que ma m\`{e}re ne la
trompait pas et que Bourbaki \'{e}tait introuvable.

\vspace{0.1cm}

\guillemotleft {\tiny \ }Pourtant elle ne renon\c{c}a pas \`{a} le d\'{e}couvrir. Elle demanda \`{a}
ses connaissances \`{a} Strasbourg s'ils savaient quelque chose de Bourbaki.
Deux ou trois seulement r\'{e}pondirent qu'ils n'en avaient jamais entendu
parler. Pour la plupart, ils croyaient bien l'avoir vu. \textquotedblleft
J'ai entendu ce nom-l\`{a}, dit une employ\'{e}e de la Biblioth\`{e}que
nationale, mais je ne peux pas mettre un visage dessus.

\vspace{0.1cm}

--- Bourbaki{\tiny \ }! Je ne connais que lui, dit un jeune stagiaire en se grattant
l'oreille. Mais je ne saurais pas vous dire qui c'est.\textquotedblright\ Le
renseignement le plus pr\'{e}cis vint de monsieur Blaise, receveur de
l'enregistrement, qui d\'{e}clara avoir eu recours aux services de Bourbaki,
l'ann\'{e}e o\`{u} son fils allait \'{e}chouer au baccalaur\'{e}at.

\vspace{0.1cm}

\guillemotleft {\tiny \ }Un matin, Mme Shuster tomba en soufflant dans le bureau de mon beau-p\`{e}re
: \textquotedblleft Je viens de voir Bourbaki. --- Ah{\tiny \ }! --- Je l'ai vu. ---
Vous croyez{\tiny \ }? --- J'en suis s\^{u}re. Il descendait du tramway, Place Kl\'{e}%
ber. Puis il a continu\'{e} \`{a} pied vers la rue des Francs-Bourgeois, il
marchait vite. Je l'ai perdu. --- \'{E}tait-ce bien lui? --- Sans aucun
doute. Un homme d'une cinquantaine d'ann\'{e}es, maigre, vo\^{u}t\'{e},
l'air d'un instituteur us\'{e}. --- Il est vrai, dit mon beau-p\`{e}re, que
ce signalement peut s'appliquer \`{a} Bourbaki. --- Vous voyez bien{\tiny \ }!
D'ailleurs, je l'ai appel\'{e}. J'ai cri\'{e}{\tiny \ }: \textit{Bourbaki{\tiny \ }!} et il
s'est retourn\'{e}. --- C'est le moyen, dit mon beau-p\`{e}re, que les
agents de la S\^{u}ret\'{e} emploient pour s'assurer de l'identit\'{e} des
malfaiteurs qu'ils recherchent. --- Quand je vous le disais, que c'\'{e}tait
lui{\tiny \ }!... J'ai bien su le trouver, moi, votre Bourbaki. Eh bien, c'est un
homme de mauvaise mine. Vous avez \'{e}t\'{e} bien imprudents, vous et votre
femme, de l'accueillir chez vous. Je me connais en physionomies et quoique
je ne l'aie vu que de dos, je jurerais que c'est une mauvaise personne. Ses
oreilles ne sont point ourl\'{e}es, et c'est un signe qui ne trompe point.
--- Ah{\tiny \ }! vous avez remarqu\'{e} que ses oreilles n'\'{e}taient point ourl%
\'{e}es{\tiny \ }? --- Rien ne m'\'{e}chappe. Mon cher Andr\'{e}, si vous ne voulez
point vous attirer d'ennuis avec votre petite famille, ne laissez plus
entrer Bourbaki chez vous, et surtout, qu'il ne vous rende point visite dans
votre lieu de travail.\textquotedblright 

\vspace{0.1cm}

\guillemotleft {\tiny \ }Or, \`{a} quelques jours de l\`{a}, il advint \`{a} M. Wickersheimer,
administrateur de la Biblioth\`{e}que nationale et grand ami des Shuster,
que les trois premiers volumes de l'\OE uvre de Gauss disparaissent des
rayons. Ayant appris la nouvelle, Mme Shuster soup\c{c}onna de suite
Bourbaki. On avait not\'{e} par le pass\'{e} la disparition de pages
illustratives dans des ouvrages d'histoire, de botanique ou d'anatomie, mais
cette fois, le vol semblait commis par un professionnel, et avec une adresse
singuli\`{e}re. Comme c'est une \oe uvre math\'{e}matique, ce ne pouvait 
\^{e}tre que lui, Bourbaki. Mme Shuster \'{e}tait devenue si obs\'{e}d\'{e}e
par notre personnage qu'elle se proposa de commander des exemplaires aupr%
\`{e}s de la Soci\'{e}t\'{e} royale des sciences de G\"{o}ttingen et d'en
faire don \`{a} la Biblioth\`{e}que, pourvu qu'on mette la main sur cet
oiseau-l\`{a}.

\vspace{0.1cm}

\guillemotleft {\tiny \ }Les Derni\`{e}res Nouvelles d'Alsace consacra m\^{e}me un article \`{a}
cette affaire o\`{u} il fut l'occasion de rendre un petit hommage au prince
des math\'{e}maticiens. Curieusement, les lecteurs s'int\'{e}ress\`{e}rent
davantage au portrait de Bourbaki qu'\`{a} celui de Gauss. Bourbaki fait d%
\'{e}sormais de l'ombre \`{a} Gauss{\tiny \ }! Et d'apr\`{e}s les renseignements
fournis en ville, \textquotedblleft il a, disait le journal, le front bas,
les yeux vairons, le regard fuyant, une patte d'oie \`{a} la tempe, les
pommettes aigu\"{e}s, rouges et luisantes. Les oreilles ne sont point ourl%
\'{e}es. Maigre un peu vo\^{u}t\'{e}, d\'{e}bile en
apparence...\textquotedblright\ 

\bigskip

\bigskip

\begin{center}
III

\bigskip

\bigskip
\end{center}

Vers dix heures du soir, Margareth ayant regagn\'{e} la chambre qu'on lui a
am\'{e}nag\'{e}e pour l'occasion, Mlle Weil dit \`{a} son fr\`{e}re{\tiny \ }:

\guillemotleft {\tiny \ }N'oublie pas de raconter comment Bourbaki s\'{e}duisit l'employ\'{e}e de la
Biblioth\`{e}que nationale.

\vspace{0.1cm}

--- J'y songeais, ma s\oe ur, r\'{e}pondit M. de Possel. L'omettre serait
perdre le plus beau de l'histoire. Elle se nommait Gudule... mais
voudriez-vous bien mesdemoiselles nous \'{e}pargner cette histoire pour le
moment{\tiny \ }? sollicita M. de Possel, en ne pouvant se retenir de rire. Revenons 
\`{a} notre Bourbaki qui, soigneusement recherch\'{e} par la justice mais
toujours introuvable, invitait chacun \`{a} mettre son amour-propre \`{a} le
trouver{\tiny \ }; les gens malin y r\'{e}ussirent. Et comme il y avait beaucoup de
gens malins \`{a} Strasbourg et aux environs, Bourbaki \'{e}tait vu en m\^{e}%
me temps \`{a} la Petite France, \`{a} la Neustadt et dans le Vieux
Cronenbourg. Un trait fut ainsi ajout\'{e} \`{a} son caract\`{e}re. On lui
accorda ce don d'ubiquit\'{e} que poss\`{e}dent tant de h\'{e}ros
populaires. Un \^{e}tre capable de franchir en un moment de longues
distances, et qui se montre tout \`{a} coup \`{a} l'endroit o\`{u} on
l'attendait le moins, effraye justement.

\vspace{0.1cm}

\guillemotleft {\tiny \ }Ainsi r\'{e}pandu dans la cit\'{e} et les environs, il restait attach\'{e} 
\`{a} notre famille par mille liens subtils. Il passait devant notre porte
et l'on croit qu'il montait parfois les escaliers de notre immeuble. On ne
le voyait jamais en face. Mais \`{a} tout moment nous connaissions son
ombre, sa voix, les traces de ses pas. Plus d'une fois nous cr\^{u}mes voir
son dos dans le cr\'{e}puscule, au tournant d'un chemin. Avec moi, il
changeait un peu de caract\`{e}re et devenait pu\'{e}ril et tr\`{e}s na\"{\i}%
f. Il se faisait moins r\'{e}el et, j'ose dire, plus po\'{e}tique. Il
entrait dans le cycle ing\'{e}nu des traditions enfantines. Il tournait au
Croquemitaine, au p\`{e}re Fouettard et au marchant de sable qui ferme, le
soir, les yeux des petits enfants. Ce n'\'{e}tait pas ce lutin qui emm\^{e}%
le, la nuit, dans l'\'{e}curie la queue des poulains. Moins rustique et
moins charmant, mais \'{e}galement espi\`{e}gle avec candeur, il fera plus
tard des moustaches d'encre \`{a} vos poup\'{e}es. Dans mon lit, avant de
m'endormir, je l'\'{e}coutais : il pleurait sur les toits avec les chats, il
aboyait avec les chiens, il emplissait de g\'{e}missements les tr\'{e}mies
et imitait dans la rue les chants des ivrognes attard\'{e}s.

\vspace{0.1cm}

\guillemotleft {\tiny \ }Ce qui me rendait Bourbaki pr\'{e}sent et familier, ce qui m'int\'{e}ressait 
\`{a} lui, c'est que son souvenir \'{e}tait associ\'{e} \`{a} tous les
objets qui m'entouraient. Mes cahiers d'\'{e}colier, dont il avait tant de
fois embrouill\'{e} et barbouill\'{e} les pages - n'\'{e}tait-ce pas le
tuteur un peu maladroit qui s'occupait de moi{\tiny \ }? Et puis le mur du balcon au-dessus
duquel j'avais vu luire, dans l'ombre, ses yeux rouges, le pot de fa\"{\i}%
ence bleue qu'une nuit d'hiver il avait fendu, \`{a} moins que ce ne f\^{u}t
la gel\'{e}e{\tiny \ }; les arbres, les rues, les bancs, tout me rappelait
Bourbaki, mon Bourbaki, le Bourbaki des enfants, des \'{e}coliers, \^{e}tre
local et mythique. Il n'\'{e}galait pas en gr\^{a}ce et en po\'{e}sie le
plus lourd \'{e}gipan, le faune le plus \'{e}pais de Sicile ou de Thessalie.
Mais c'\'{e}tait un demi-dieu encore.

\vspace{0.1cm}

\guillemotleft {\tiny \ }Pour mon beau-p\`{e}re, il avait un tout autre caract\`{e}re : il \'{e}tait
embl\'{e}matique et philosophique. M. Weil avait une grande piti\'{e} des
hommes. Il ne les croyait pas tr\`{e}s raisonnables{\tiny \ }; leurs
erreurs, quand elles n'\'{e}taient point cruelles, l'amusaient et le
faisaient sourire. La croyance en Bourbaki l'int\'{e}ressait comme un abr%
\'{e}g\'{e} et un compendium de toutes les croyances humaines. Comme il \'{e}%
tait ironique et moqueur, il parlait de Bourbaki ainsi que d'un \^{e}tre r%
\'{e}el. Il y mettait parfois tant d'insistance et marquait les
circonstances avec une telle exactitude, que ma m\`{e}re en \'{e}tait toute
surprise et lui disait dans sa candeur : \textquotedblleft On dirait que tu
parles s\'{e}rieusement, mon ami{\tiny \ }: tu sais pourtant
bien...\textquotedblright 

\vspace{0.1cm}

\guillemotleft {\tiny \ }Il r\'{e}pondait gravement{\tiny \ }: \textquotedblleft Tout le monde croit \`{a}
l'existence de Bourbaki. Serais-je un bon citoyen si je la niais{\tiny \ }? Il faut y
regarder \`{a} deux fois avant de supprimer un article de la foi
commune.\textquotedblright 

\vspace{0.1cm}

\guillemotleft {\tiny \ }Un esprit parfaitement honn\^{e}te a seul de semblables scrupules. Au fond,
M. Weil \'{e}tait gassendiste. Il accordait son sentiment particulier avec
le sentiment public, croyant comme les Alsaciens \`{a} l'existence de
Bourbaki, mais n'admettant pas son intervention directe dans le vol des
livres et la s\'{e}duction des femmes. Enfin il professait sa croyance en
l'existence d'un Bourbaki, pour \^{e}tre bon Alsacien{\tiny \ }; et il se
passait de Bourbaki pour expliquer les \'{e}v\'{e}nements qui
s'accomplissaient dans la ville. De sorte qu'en cette circonstance, comme en
toute autre, il fut un galant homme et un bon esprit.

\vspace{0.1cm}

\guillemotleft {\tiny \ }Quant \`{a} notre m\`{e}re, elle se reprochait un peu la naissance de
Bourbaki, et non sans raison. Car enfin Bourbaki \'{e}tait n\'{e} d'un
mensonge de notre m\`{e}re, comme Caliban du mensonge du po\`{e}te. Sans
doute les fautes n'\'{e}taient pas \'{e}gales et ma m\`{e}re \'{e}tait plus
innocente que Shakespeare. Pourtant elle \'{e}tait effray\'{e}e et confuse
de voir son mensonge bien mince grandir d\'{e}mesur\'{e}ment, et sa l\'{e}g%
\`{e}re imposture remporter un si prodigieux succ\`{e}s, qui ne s'arr\^{e}%
tait pas, qui s'\'{e}tendait sur toute une ville et mena\c{c}ait de s'\'{e}%
tendre sur le monde. Un jour m\^{e}me elle p\^{a}lit, croyant qu'elle allait
voir son mensonge se dresser devant elle. Ce jour-l\`{a}, une bonne qu'elle
avait, nouvelle dans la maison et dans le pays, vint lui dire qu'un homme
demandait \`{a} la voir. Il avait, disait-il, besoin de parler \`{a} madame.
\textquotedblleft Quel homme est-ce{\tiny \ }? --- Un homme en costume crois\'{e},
avec un cartable en cuir. Il a l'air d'un professeur. --- A-t-il dit son nom{\tiny \ }? --- Oui, madame. --- Eh bien, comment se nomme-t-il{\tiny \ }? --- Bourbaki. --- Il
vous a dit qu'il se nommait{\tiny \ }?... --- Bourbaki, oui, madame. --- Il est ici{\tiny \ }?... --- Oui, madame. Il attend dans le salon. --- Vous l'avez vu{\tiny \ }? --- Oui,
madame. --- Qu'est-ce qu'il veut{\tiny \ }? --- Il ne me l'a pas dit. Il ne veut le
dire qu'\`{a} madame. --- Allez le lui demander.\textquotedblright 

\vspace{0.1cm}

\guillemotleft {\tiny \ }Quand la servante retourna dans le salon, Bourbaki n'y \'{e}tait plus. Cette
rencontre de la servante \'{e}trang\`{e}re et de Bourbaki ne fut jamais \'{e}%
claircie. Mais je crois qu'\`{a} partir de ce jour ma m\`{e}re commen\c{c}a 
\`{a} croire que Bourbaki pouvait bien exister, et que les Weil pouvaient
bien n'avoir pas menti.{\tiny \ }\guillemotright 

\bigskip

\bigskip

\bigskip

\begin{center}
* * *

\bigskip

\newpage
\end{center}

\textit{On l'aura donc reconnu, il s'agit d'un d\'{e}tournement, au bon sens
du terme s'il convient d'en inventer un, d'un c\'{e}l\`{e}bre conte
satirique d'Anatole France{\tiny \ }\footnote{%
\noindent \textit{Anatole-Fran\c{c}ois Thibault (Anatole France),
\textquotedblleft Crainquebille, Putois, Riquet et plusieurs autres r\'{e}%
cits profitables\textquotedblright , Calmann-L\'{e}vy, 1904.}}, o\`{u} le
personnage de Putois, n\'{e} d'un mensonge de Mme Bergeret, a d\^{u} laisser
place \`a Nicolas Bourbaki, dont nous connaissons tous l'histoire et ses anecdotes. On l'aura aussi remarqu\'{e}, la date correspond au centenaire
de la mort du grand homme de lettres que fut Anatole France (16 avril 1844 -
12 octobre 1924). Ce pastiche historico-math\'{e}matico-litt\'{e}raire met \`{a} l'honneur la famille Weil en les personnes d'Andr\'{e}
Weil (\'{E}loi Bergeret), sa femme \'{E}velyne (Mme Bergeret) et leurs deux
filles, Sylvie (Mlle Zo\'{e} Bergeret) et Nicolette (jouant le r\^{o}le de
MM Goubin et Jean Marteau, permettant ainsi \`{a} la conversation de se d%
\'{e}rouler dans une ambiance familiale). Est aussi \`{a} l'honneur la
famille de{\tiny \ }Possel-Deydier avec Alain de Possel (Lucien Bergeret), fils de Ren%
\'{e} de Possel et \'{E}velyne Gillet, avant que Mme de Possel ne devienne
Mme Weil, et sa fille cadette Margareth (Pauline) qu'il a eue avec \'{E}lise
Lucciani, \'{e}pous\'{e}e en secondes noces. Les lieux ont \'{e}t\'{e} aussi
chang\'{e}s{\tiny \ }; Saint-Omer devient Strasbourg et, \`{a} quelques
lieues, Monplaisir devient Brumath, Andr\'{e} Weil ayant exerc\'{e} \`{a}
l'universit\'{e} de Strasbourg, avec son ami Henri Cartan, durant la p\'{e}%
riode qui a vu na\^{\i}tre Bourbaki. Pour rester fid\`{e}le au texte
original d'Anatole France, il fallait deux personnages cl\'{e}s{\tiny \ }: le fr\`{e}%
re et la s\oe ur qui parlent de leurs parents. Nous avons donc suppos\'{e}
que le petit Alain a v\'{e}cu avec sa m\`{e}re et son beau-p\`{e}re \`{a}
Strasbourg, alors qu'il se pourrait bien qu'il ait v\'{e}cu avec son p\`{e}%
re et sa belle-m\`{e}re, Yvonne Lib\'{e}rati. Les Souvenirs{\tiny \ }\footnote{%
\noindent \textit{Andr\'{e} Weil, \textquotedblleft Souvenirs
d'apprentissage\textquotedblright , Birkh\"{a}user, 1991.}} viennent nous
confirmer, justement, qu'apr\`{e}s son retour des \'{E}tats-Unis, \`{a} la
rentr\'{e}e universitaire 1937, Andr\'{e} Weil \textquotedblleft \lbrack
devait] trouver un logis \`{a} Strasbourg pour \'{E}velyne, son fils Alain
alors \^{a}g\'{e} de six ans, et [lui]\textquotedblright . Ainsi, en \'{e}%
voquant le personnage de Bourbaki, tout au d\'{e}but de l'histoire, Alain
(Lucien), qui \'{e}tait t\'{e}moin de la naissance de Bourbaki le tuteur
(Putois, le jardinier), s'adressait \`{a} sa demi-s\oe ur Sylvie (Zo\'{e}%
), ainsi qu'\`{a} sa fille Margareth (Pauline). Pouvait-on mieux faire{\tiny \ }? La
description un peu pittoresque de la situation de l'appartement des Weil,
tout au d\'{e}but du texte, est emprunt\'{e}e \`{a} Pierre Cartier (qui nous
a quitt\'{e}s r\'{e}cemment), dans son vibrant hommage \`{a} Andr\'{e} Weil{\tiny \ }\footnote{%
\noindent \textit{Pierre Cartier, \textquotedblleft Andr\'{e} Weil
(1906-1998){\tiny \ }: adieu \`{a} un ami\textquotedblright , S\'{e}minaire de
Philosophie et Math\'{e}matiques, numdam.org, 1998.}}. On sait qu'il est situ%
\'{e} rue Auguste-Comte, dans le quartier de l'Od\'{e}on du 6}$^{\mathit{e}}$%
\textit{\ arrondissement de Paris. Mais pour une tout autre immersion dans
le monde des Weil, en particulier dans \textquotedblleft le g\'{e}nie bic%
\'{e}phale\textquotedblright\ d'Andr\'{e} et de Simone, une autre adresse bien
connue{\tiny \ }: Chez les Weil{\tiny \ }\footnote{%
\noindent \textit{Sylvie Weil, \textquotedblleft Chez les Weil{\tiny \ }: Andr\'{e}
et Simone\textquotedblright , Buchet-Chastel, 2009.}}, dont l'illustration
de couverture, simple et innocente, nous rappelle un peu la l\'{e}gendaire
collection \textquotedblleft Lis{\tiny \ }!\textquotedblright ,\ de notre cher et
regrett\'{e} Ma\^{\i}tre Boukmakh{\tiny \ }\footnote{\textit{Il s'agit en fait de la collection \textquotedblleft Iqrae\textquotedblright , qui a marqué l'imaginaire de plusieurs g\'{e}n\'{e}rations d'\'{e}l\`{e}ves au Maroc. Le point d'exclamation,
ajout\'{e}, \textquotedblleft est une invitation \`{a} lire. C'est un point
d'excitation, d'incitation, d'enthousiasme, mais pas
d'ordre\textquotedblright , diraient Bernard Pivot et sa fille C\'{e}cile{\tiny \ }! }%
}{\tiny \ }: le fr\`{e}re (tenant un livre) et sa s\oe ur, assis c\^{o}te \`{a} c\^{o}%
te, comme deux \'{e}coliers en classe. }

\bigskip

\textit{Comme on peut le constater, certains 
\'{e}v\'{e}nements de l'histoire originale concernant Putois, le p\'{e}pini%
\'{e}riste, ont \'{e}t\'{e} transform\'{e}s pour mieux correspondre au
personnage de Bourbaki, le math\'{e}maticien. Par exemple, Mme Cornouiller
qui attendait d\'{e}sesp\'{e}r\'{e}ment Putois pour tailler son jardin dans
son manoir de Monplaisir est devenue Mme Schuster qui cherchait Bourbaki
pour le faire travailler dans une \'{e}cole \`{a} Brumath. Aussi, les trois
melons de Mme Cornouiller ont \'{e}t\'{e} remplac\'{e}s par les trois
premiers volumes de l'\OE uvre de Gauss dans la Biblioth\`{e}que nationale
de Strasbourg{\tiny \ }(!) Si le personnage de Mme Schuster est fictif, celui de M.
Wickersheimer, mis \`{a} l'honneur, est bien r\'{e}el car c'est lui qui
occupait la fonction d'administrateur de la Biblioth\`{e}que de 1926 \`{a}
1950 (et \`{a} Clermont-Ferrand de 1941 \`{a} 1945). Quant \`{a} l'affaire
Gudule, entre autres faits divers, si l'omettre serait perdre le plus beau
de l'histoire, comme le faisait dire Anatole France \`{a} Lucien Bergeret,
nous avons tout simplement estim\'{e} qu'elle ne correspondrait pas au
personnage de Bourbaki. }

\bigskip

\textit{Enfin, par-del\`{a} la question de savoir lequel des Nicolas a bel
et bien exist\'{e}, le Nicolas de Sylvie Vartan, celui de William Sheller ou
encore notre Nicolas \`{a} nous, Nicolas Bourbaki, le pr\'{e}sent texte se
veut d'une part un hommage \`{a} Anatole France, cent ans apr\`{e}s sa mort
(et presque autant, un peu plus, apr\`{e}s le canular de Raoul Husson \`{a}
la Normale Sup'), et d'autre part, un hommage \`{a} Nicolas Bourbaki et \`{a}
ses \textquotedblleft Godparents\textquotedblright , Andr\'{e} et \'{E}%
velyne Weil. Nous esp\'{e}rons ne pas avoir franchi les limites de l'humour
courtois et bienveillant, tenant en haute estime toutes les personnes dont
les noms, ou ceux de leurs proches, ont \'{e}t\'{e} \'{e}voqu\'{e}s, y compris, bien entendu, la famille Bourbaki{\tiny \ }! Un
dernier brin de fantaisie pour finir{\tiny \ }: dans le parcours scientifique tout
passionnant d'Andr\'{e} Weil, depuis l'Europe jusqu'\`{a} Princeton, en
passant par l'Inde et le Br\'{e}sil, on sait que lorsque les Weil, avec le
petit Alain, ont embarqu\'{e} pour Fort-de-France \`{a} bord du Winnipeg,
pour rejoindre ensuite New-York, ils ont d\^{u} faire une escale de trois
jours \`{a} Casablanca. On ne peut d\`{e}s lors s'emp\^{e}cher d'imaginer M.
Weil, attabl\'{e} au Rick's Caf\'{e} Am\'{e}ricain, berc\'{e} par
\textquotedblleft As Time Goes By\textquotedblright , sous le regard
suspicieux du capitaine Louis Renault, en train de feuilleter son
\textquotedblleft Int\'{e}gration dans les groupes
topologiques\textquotedblright\ qui vient de para\^{\i}tre, puis chantant La
Marseillaise avec Victor Laszlo, au grand dam du major Heinrich Strasser{\tiny \ }! }

\bigskip

\bigskip

\ \ \ \ \ \ \ \ \ \ \ \ \ \ \ \ \ \ \ \ \ \ \ \ \ \ \ \ \ \ \ \ \ \ \ \ \ \
\ \ \ \ \ \ \ \ \ \ \ \ \ \ \ \ \textit{\ \ \ \ \ \ \ \ \ \ \ \ \ \ \ \ \ \
\ \ \ \ \ Rabat, automne 2024}

\end{document}